\documentclass[10pt]{article}
\usepackage[psamsfonts]{amssymb}
\usepackage{amsthm,amsmath}
\usepackage[all]{xy}

\title{\bf Meaningful aggregation functions mapping ordinal scales into an ordinal scale: a state of the art}

\author{%
Jean-Luc MARICHAL \\
Mathematics Research Unit, University of Luxembourg \\
162A, avenue de la Fa\"{\i}encerie, L-1511 Luxembourg, Luxembourg\\
jean-luc.marichal[at]uni.lu%
\and %
Radko MESIAR \\
Department of Mathematics, Slovak University of Technology \\
Radlinskeho 11, SK-81368 Bratislava, Slovakia \\
mesiar[at]math.sk%
}

\date{Revised, March 13, 2009}

\begin{document}
\maketitle

\theoremstyle{plain}
\newtheorem{theorem}{Theorem}[section]
\newtheorem{lemma}{Lemma}[section]
\newtheorem{proposition}{Proposition}[section]
\newtheorem{corollary}{Corollary}[section]

\theoremstyle{definition}
\newtheorem{definition}{Definition}[section]
\newtheorem{example}{Example}[section]
\newtheorem{conjecture}{Conjecture}[section]
\newtheorem{remark}{Remark}[section]

\newcommand{\N}{\mathbb{N}}
\newcommand{\R}{\mathbb{R}}
\newcommand{\Vspace}{\vspace{2ex}}
\newcommand{\ran}{{\rm ran}}
\newcommand{\ien}{\mathcal{I}_n[E]}
\newcommand{\isen}{\mathcal{I}_n^*[E]}
\newcommand{\lc}{L_c}
\newcommand{\lxi}{L_{\xi(I)}}
\newcommand{\cm}{\mathcal{C}_\mu}
\newcommand{\sm}{\mathcal{S}_\mu}
\newcommand{\cn}{\mathcal{C}_n}
\newcommand{\cen}{\mathcal{C}_n[E]}
\newcommand{\bx}{\mathbf{x}}
\newcommand{\by}{\mathbf{y}}
\newcommand{\ba}{\mathbf{a}}
\newcommand{\bb}{\mathbf{b}}
\newcommand{\bff}{\mathbf{f}}
\newcommand{\bphi}{\boldsymbol{\phi}}
\newcommand{\loe}{~\Big\{{<\atop =}\Big\}~}
\def\argmax{\mathop{\rm argmax}}
\def\bigtimes{\mathop{\boldsymbol{\times}}}

\begin{abstract}
We present an overview of the meaningful aggregation functions mapping ordinal scales into an ordinal scale. Three main classes are discussed,
namely order invariant functions, comparison meaningful functions on a single ordinal scale, and comparison meaningful functions on independent
ordinal scales. It appears that the most prominent meaningful aggregation functions are lattice polynomial functions, that is, functions built
only on projections and minimum and maximum operations.
\end{abstract}


%---------------------------------------------------------------------------------------------- Section 1
\section{Introduction}

In many domains we are faced with the problem of aggregating a collection of numerical readings to obtain an average value. Actually, such an
aggregation problem is becoming more and more present in an increasing number of areas not only of mathematics or physics, but also of
engineering, economical, social, and other sciences. Various aggregation functions and processes have already been proposed in the literature
and many others are still to be designed to fulfill newer and newer requirements.

Studies on the aggregation problem have shown that the choice of the aggregation function is far from being arbitrary and should be based upon
properties dictated by the framework in which the aggregation is performed. One of the main concerns when choosing an appropriate function is to
take into account the scale types of the variables being aggregated. On this issue, Luce~\cite{Luc59} observed that the general form of the
functional relationship between variables is greatly restricted if we know the scale types of the dependent and independent variables. For
instance, if all the variables define a common ordinal scale, it is clear that any relevant aggregation function cannot be constructed from
usual arithmetic operations, unless these operations involve only order. Thus, computing the arithmetic mean is forbidden whereas the median or
any order statistic is permitted.

Specifically, suppose $x_1,\ldots,x_n,x_{n+1}$ are $n+1$ variables, each $x_i$ having a real interval as a domain, and $x_{n+1}=F(x_1,\ldots,x_n)$ is some unknown aggregation function. The problem is to find the general form of the
function $F$ knowing the scale types of the input and output variables. The \emph{scale type}\/ of a variable $x_i$ is defined by the class of
\emph{admissible transformations}, transformations, such as that from grams to pounds or degrees Fahrenheit to degrees centigrade, that change
the scale into an alternative acceptable scale. In the case of a \emph{ratio scale}, for example, an admissible transformation is a function of
the form $x\mapsto rx$, with some $r>0$, which changes the unit of the scale. Similarly, for an \emph{interval scale}, an admissible
transformation is a function $x\mapsto rx+s$, with $r>0$ and $s\in\R$, which modifies both the origin and the unit of the scale. For an
\emph{ordinal scale}, an admissible transformation is a strictly increasing function $x\mapsto \phi(x)$, which changes the values of the scale
while preserving their order. For more details on the theory of scale types, see \cite{KraLucSupTve71,LucKraSupTve90,Nar02,Rob79,Rob94}.

Luce's principle, called ``principle of theory construction", is based on the requirement that admissible transformations of the input variables
must lead to an admissible transformation of the output variable. For example, if the input variables are independent scales, then the
aggregation function $F$ should satisfy the following condition. For any admissible transformations $\phi_1,\ldots,\phi_n$ of the input
variables, there is an admissible transformation $\psi_{\phi_1,\ldots,\phi_n}$ of the output variable so that $F\big(\phi_1(x_1),\ldots,\phi_n(x_n)\big)=\psi_{\phi_1,\ldots,\phi_n}(x_{n+1})$ or, equivalently,
\begin{equation}\label{eq:fppf1}
F\big(\phi_1(x_1),\ldots,\phi_n(x_n)\big)=\psi_{\phi_1,\ldots,\phi_n}\big(F(x_1,\ldots,x_n)\big).
\end{equation}

The solutions of this functional equation constitute the set of the possible aggregation functions, which are ``meaningful" in the sense that
they do not depend upon the particular scales of measurement chosen for the variables, but only upon their scale types.

We can also assume that the input variables define the same scale, which implies that the same admissible transformation must be applied to all
the input variables. In this case, the condition on $F$ is the following. For any common admissible transformation $\phi$ of the input
variables, there is an admissible transformation $\psi_{\phi}$ of the output variable so that
\begin{equation}\label{eq:fppf2}
F\big(\phi(x_1),\ldots,\phi(x_n)\big)=\psi_{\phi}\big(F(x_1,\ldots,x_n)\big).
\end{equation}

In the extreme case where all the input and output variables define the same scale, then, for any admissible transformation $\phi$ of the input
and output variables, we must have
\begin{equation}\label{eq:fppf3}
F\big(\phi(x_1),\ldots,\phi(x_n)\big)=\phi\big(F(x_1,\ldots,x_n)\big).
\end{equation}

Equations (\ref{eq:fppf1}) and (\ref{eq:fppf2}) were completely solved in the eighties for ratio scale variables, with domain
$\R_+:=[0,\infty)$, and interval scale variables, with domain $\R$, even under some further assumptions such as symmetry, continuity, and
nondecreasing monotonicity (in each argument); see \cite{Acz87,AczGroSch94,AczRob89,AczRobRos86,Kim90}.

\begin{example}[\cite{AczRobRos86}]
If all the input variables are independent ratio scales and the output variable is also a ratio scale then the meaningful aggregation functions
$F\colon\R^n_+\to\R_+$ are exactly the solutions of (\ref{eq:fppf1}), where each admissible transformation is a multiplication by a positive
constant. These solutions are given by
$$F(x_1,\ldots,x_n)=a\prod_{i=1}^n f_i(x_i)\qquad (a>0),$$ where the functions $f_i\colon\R_+\to\R_+$ fulfill the
equations $f_i(x_iy_i)=f_i(x_i)f_i(y_i)$ $(i=1,\ldots,n)$. Under continuity, these solutions are of the form
$$
F(x_1,\ldots,x_n)=a\prod_{i=1}^n x_i^{c_i}\qquad (a>0, \; c_1,\ldots,c_n\in\R).
$$
\end{example}

For ordinal scales and without any further assumptions, equations (\ref{eq:fppf1}) and (\ref{eq:fppf2}) have resisted and have remained unsolved
for a long time. At present, the complete description of their solutions is known and has been presented recently in a couple of articles; see \cite{MarMes04,MarMesRuc05}.

The main purpose of this paper is to present a catalog of all the possible meaningful aggregation functions mapping ordinal scales into an
ordinal scale. More precisely, we yield all the possible solutions of each of the functional equations above, where the admissible
transformations are strictly increasing functions. We also present the solutions under some further assumptions such as continuity, symmetry,
idempotency, and nondecreasing monotonicity (in each argument).

In such an ordinal framework, it is natural to assume that the common domain of the input variables be any open real interval or even the whole
real line. However, we consider the more general situation where the domain of the input variables is any real interval $E$, possibly unbounded,
and where the domain of the output variable is the real line, except for equation (\ref{eq:fppf3}) where this domain must also be the set $E$.
We further assume that the admissible transformations of the input variables are confined to the increasing bijections from $E$ onto $E$. This
latter assumption, which brings no restriction to the solutions of the functional equations, enables us to consider closed domains $E$ whose
endpoints remain fixed under any admissible transformation.

We thus provide and discuss all the solutions $F\colon E^n\to\R$ of the functional equations (\ref{eq:fppf1}) and (\ref{eq:fppf2}), where
$\phi_1,\ldots,\phi_n$, and $\phi$ are arbitrary increasing bijections from $E$ onto $E$ and where $\psi_{\phi_1,\ldots,\phi_n}$ and
$\psi_{\phi}$ are strictly increasing functions from $\ran(F)$ into $\ran(F)$. We call those solutions \emph{comparison meaningful functions on
independent ordinal scales}\/ and \emph{comparison meaningful functions on a single ordinal scale}, respectively. We also provide all the
solutions $F\colon E^n\to E$ of the functional equation (\ref{eq:fppf3}), where $\phi$ is an arbitrary increasing bijection from $E$ onto $E$. We call
those solutions \emph{order invariant functions}.

The outline of this paper is as follows. In \S{2} we introduce the concept of order invariant subsets, which will play a key role in the
description of the most general solutions of the three functional equations. In \S{3} we introduce the lattice polynomial functions, which
represent most of the regular (e.g., nondecreasing) solutions of equations (\ref{eq:fppf2}) and (\ref{eq:fppf3}). We also recall some of their
properties as aggregation functions. In \S{4} we present and discuss all the order invariant functions. In \S\S{5} and 6 we present respectively
the comparison meaningful functions on a single ordinal scale and the comparison meaningful functions on independent ordinal scales. Finally in
\S{7} we provide interpretations of equations (\ref{eq:fppf1})--(\ref{eq:fppf3}) in the setting of aggregation on finite chains.

Throughout this paper we denote by $E$ any real interval, bounded or not, with interior $E^\circ$. We also denote by $B[E]$ the set of
\emph{included boundaries} of $E$, that is $B[E]:=E\setminus E^\circ$. The set of all increasing bijections $\phi$ of $E$ onto itself is denoted
by $\Phi[E]$. As each function $\phi\in\Phi[E]$ preserves the ordinal structure of $E$, the set $\Phi[E]$ is actually the \emph{order
automorphism group}, under composition, of $E$. Finally, the symbol $[n]$ denotes the index set $\{1,\ldots,n\}$ and, for any $\bx\in E^n$ and
any $\bphi\in\Phi[E]^n$, the symbol $\bphi(\bx)$ denotes the vector $\big(\phi_1(x_1),\ldots,\phi_n(x_n)\big)$.

%---------------------------------------------------------------------------------------------- Section 2
\section{Order invariant subsets}\label{sec:gfass3}

The space $E^n$ can be partitioned into \emph{order invariant subsets}, which are very useful in describing the general solutions of the
functional equations introduced above. Those subsets were introduced first by Ovchinnikov (see \cite[\S{3}]{Ovc96} and \cite[\S{2}]{Ovc98c}) in
the general framework of ordered sets and then independently by Bart{\l}omiejczyk and Drewniak~\cite{BarDre04} for closed real intervals; see
also \cite{MarMes04,MarMesRuc05,MesRuc04}. In this section we introduce them through the concept of group orbit.\footnote{For definitions and
results about the concept of orbit in algebra, see e.g.\ \cite{Gri07}.}

Consider the product set $\Phi[E]^n$ and its \emph{diagonal restriction}
$$\Phi_n[E]:=\big\{(\underbrace{\phi,\ldots,\phi}_n) : \phi\in \Phi[E]\big\}.$$ As $\Phi_n[E]$ is clearly a subgroup of $\Phi[E]^n$, we can
define the orbit of any element $\bx\in E^n$ under the action of $\Phi_n[E]$, that is, $\Phi_n[E](\bx) := \{\bphi(\bx) : \bphi \in \Phi_n[E]\}$. The set of orbits of $E^n$ under $\Phi_n[E]$ forms a partition of $E^n$ into equivalence classes, where $\bx,\by\in E^n$ are equivalent if their
orbits are the same, that is, if there exists $\bphi\in\Phi_n[E]$ such that $\by=\bphi(\bx)$. The orbits of $E^n$ under $\Phi_n[E]$ are \emph{order invariant subsets} in the following sense (see \cite{BarDre04}).

\begin{definition}
A nonempty subset $I$ of $E^n$ is called \emph{order invariant}\/ if, for any $\bx\in I$, we have $\bphi(\bx)\in I$ for all $\bphi\in
\Phi_n[E]$.\footnote{Equivalently, $I$ is order invariant if $\bphi(I)\subseteq I$ for all $\bphi\in\Phi_n[E]$. Actually, since $\Phi_n[E]$ is a
group, we can even write $\bphi(I)=I$.} An order invariant subset of $E^n$ is \emph{minimal}\/ if it has no proper order invariant subset.
\end{definition}

It is easy to see that the set $\ien:=E^n/\Phi_n[E]$ of orbits of $E^n$ under $\Phi_n[E]$ is identical to the set of minimal order invariant
subsets of $E^n$. Moreover, any order invariant subset is a union of those orbits.

The following proposition (for closed $E$, see \cite{BarDre04,MesRuc04}) yields a complete description of the orbits.

\begin{proposition}\label{prop:DescOrbits}
We have $I\in \ien$ if and only if there exists a permutation $\pi$ on $[n]$ and a sequence $\{\lhd_i\}_{i=0}^n$ of symbols $\lhd_i\in\{<,=\}$,
containing at least one symbol $<$ if $\inf E\in E$ and $\sup E\in E$, such that
$$I=\{\bx\in E^n : \inf E \,\lhd_0\, x_{\pi(1)} \,\lhd_1\,\cdots\,\lhd_{n-1}\, x_{\pi(n)} \,\lhd_n\, \sup E\},$$ where
$\lhd_0$ is $<$ if $\inf E\notin E$ and $\lhd_n$ is $<$ if $\sup E\notin E$.
\end{proposition}

\begin{example}[\cite{MesRuc04}]\label{ex:fggger}
The unit square $[0,1]^2$ contains exactly eleven minimal order invariant subsets, namely the open triangles $\{(x_1,x_2) : 0<x_1<x_2<1\}$ and
$\{(x_1,x_2) : 0<x_2<x_1<1\}$, the open diagonal $\{(x_1,x_2) : 0<x_1=x_2<1\}$, the four square vertices, and the four open line segments
joining neighboring vertices.
\end{example}

\begin{remark}\label{rem:perm2}
From Proposition~\ref{prop:DescOrbits} we can easily derive an alternative way to characterize the membership of given vectors $\bx,\by\in E^n$
in the same orbit. Let $\Pi_n$ be the set of permutations on $\{0,1,\ldots,n+1\}$ and, for any $\bx\in E^n$, define
$$
\Pi(\bx):=\{\pi\in\Pi_n : x_{\pi(0)}\leqslant x_{\pi(1)} \leqslant\cdots\leqslant x_{\pi(n+1)}\},
$$
where $x_0:=\inf E$ and $x_{n+1}:=\sup E$. Then, for any $\bx,\by\in E^n$, there exists $I\in\ien$ such that $\bx,\by\in I$ if and only if
$\Pi(\bx)=\Pi(\by)$.\footnote{This condition is more restrictive than \emph{comonotonicity} of vectors $\bx$ and $\by$, which simply means that
$\Pi(\bx)$ and $\Pi(\by)$ overlap; see \cite{HarLitPol52}.}
\end{remark}

Since $\Phi[E]^n$ is itself a group, we can also define the orbit of any element $\bx\in E^n$ under the action of $\Phi[E]^n$, that is, $\Phi[E]^n(\bx) := \{\bphi(\bx) : \bphi \in \Phi[E]^n\}$. Just as for the subgroup $\Phi_n[E]$, the set of orbits of $E^n$ under $\Phi[E]^n$ forms a partition of $E^n$ into equivalence classes, where
$\bx,\by\in E^n$ are equivalent if there exists $\bphi\in\Phi[E]^n$ such that $\by=\bphi(\bx)$. The orbits of $E^n$ under $\Phi[E]^n$ are \emph{strongly order invariant subsets} in the following sense (see \cite{MarMesRuc05}).

\begin{definition}
A nonempty subset $I$ of $E^n$ is called \emph{strongly order invariant}\/ if, for any $\bx\in I$, we have $\bphi(\bx)\in I$ for all $\bphi\in
\Phi[E]^n$.\footnote{Equivalently, $I$ is strongly order invariant if $\bphi(I)\subseteq I$ for all $\bphi\in\Phi[E]^n$. Once again, since
$\Phi[E]^n$ is a group, we can even write $\bphi(I)=I$.} A strongly order invariant subset of $E^n$ is \emph{minimal}\/ if it has no proper
strongly order invariant subset.
\end{definition}

The set $\isen:=E^n/\Phi[E]^n$ of orbits of $E^n$ under $\Phi[E]^n$ is identical to the set of minimal strongly order invariant subsets of
$E^n$. Moreover, any strongly order invariant subset is a union of those orbits.

The following proposition \cite{MarMesRuc05} yields a complete description of the orbits.

\begin{proposition}
We have $\isen  = \{\times_{i=1}^n I_i : I_i\in \mathcal{I}_1[E]\}= (\mathcal{I}_1[E])^n$, with cardinality $|\isen|=(1+|B[E]|)^n$.
\end{proposition}

\begin{example}[\cite{MesRuc04}]
The unit square $[0,1]^2$ contains exactly nine minimal strongly order invariant subsets, namely the open square $(0,1)^2$, the four square
vertices, and the four open line segments joining neighboring vertices.
\end{example}

Let us now show that the set $\isen$ can be described by means of the set $\ien$. For any $i\in [n]$, let ${\rm P}_i\colon E^n\to E$ be the projection
operator onto the $i$th coordinate, that is, ${\rm P}_i(\bx):=x_i$. We can easily see that, for any $I\in\ien$, we have ${\rm
P}_i(I)\in\mathcal{I}_1[E]$. Define an equivalence relation $\sim$ on $\ien$ as
$$
I\sim J \quad \Leftrightarrow \quad {\rm P}_i(I)={\rm P}_i(J) \quad (i\in [n]).
$$
Then, it is easy to see \cite{MarMesRuc05} that
$$
\isen  = \Big\{\bigcup_{\textstyle{J\in\ien\atop J\sim I}}J : I\in \ien\Big\} = \Big\{\bigtimes_{i=1}^n {\rm P}_i(I) : I\in \ien\Big\}.
$$

Now, to easily describe certain nondecreasing aggregation functions, it is useful to consider partial orders on $\ien$ and $\isen$. Starting
from the natural order $\{\inf E\} \prec E^\circ \prec \{\sup E\}$ on $\mathcal{I}_1[E]$, we can straightforwardly derive a partial order
$\preccurlyeq$ on $\ien$, namely
$$
I\preccurlyeq J\quad \Leftrightarrow \quad {\rm P}_i(I)\preccurlyeq {\rm P}_i(J) \quad (i\in [n]).
$$
The corresponding partial order on $\isen$ is defined similarly.

\begin{remark}
Consider again the set $\Pi_n$ of permutations on $\{0,1,\ldots,n+1\}$ (see Remark~\ref{rem:perm2}). For any $\bx\in E^n$, we can define
\begin{eqnarray*}
\lefteqn{\Pi^*(\bx) :=\{\pi\in\Pi_n : \pi(i)\leqslant \ell(\bx) ~\Leftrightarrow ~x_i=\inf E}\\
&& \mbox{ and }~ \pi(j)\geqslant n+1-u(\bx) ~\Leftrightarrow ~x_j=\sup E\},
\end{eqnarray*}
where $x_0:=\inf E$, $x_{n+1}:=\sup E$ and $\ell(\bx):=\{i\in [n] : x_i=\inf E\}$, $u(\bx):=\{j\in [n] : x_j=\sup E\}$. Then, for any
$\bx,\by\in E^n$, there exists $I\in\isen$ such that $\bx,\by\in I$ if and only if $\Pi^*(\bx)=\Pi^*(\by)$.
\end{remark}

%---------------------------------------------------------------------------------------------- Section 3
\section{Lattice polynomial functions and some of their properties}

As we will see in the subsequent sections, certain solutions of equations (\ref{eq:fppf2}) and (\ref{eq:fppf3}) are constructed from
\emph{lattice polynomial functions}. In this section we briefly recall the basic material about these functions. As we are concerned with
aggregation functions defined in real domains, we do not consider lattice polynomial functions on a general lattice, but simply on $\R$, which
is a particular lattice. The lattice operations $\wedge$ and $\vee$ then represent the minimum and maximum operations, respectively.

\subsection{Lattice polynomial functions}
\label{sec:sfdaag}

Let us first recall the concept of lattice polynomial function (with real variables); see e.g.\ Birkhoff~\cite[\S{II.5}]{Bir67} or
Gr\"atzer~\cite[\S{I.4}]{Grae03}.

\begin{definition}
The class of lattice polynomial functions from $\R^n$ to $\R$ is defined as follows.
\begin{enumerate}
\item[(i)] For any $k\in [n]$, the projection $\mathrm{P}_k\colon \bx\mapsto x_k$ is a lattice polynomial function from $\R^n$ to $\R$.

\item[(ii)] If $p$ and $q$ are lattice polynomial functions from $\R^n$ to $\R$, then $p\wedge q$ and $p\vee q$ are lattice polynomial functions
from $\R^n$ to $\R$.

\item[(iii)] Every lattice polynomial function from $\R^n$ to $\R$ is constructed by finitely many applications of the rules (i) and (ii).
\end{enumerate}
\end{definition}

Because $\R$ is a distributive lattice, any lattice polynomial function can be written in {\em disjunctive}\/ and {\em conjunctive}\/ forms as
follows; see e.g.\ \cite[\S{II.5}]{Bir67}.

\begin{proposition}\label{prop:lp dnf}
Let $p\colon\R^n\to \R$ be any lattice polynomial function. Then there are nonconstant set functions $\alpha\colon 2^{[n]}\to\{0,1\}$ and
$\beta\colon 2^{[n]}\to\{0,1\}$, with $\alpha(\varnothing)=0$ and $\beta(\varnothing)=1$, such that
\begin{equation}\label{eq:pab}
p(\bx)=\bigvee_{\textstyle{S\subseteq [n]\atop \alpha(S)=1}}\bigwedge_{i\in S}x_i = \bigwedge_{\textstyle{S\subseteq [n]\atop
\beta(S)=0}}\bigvee_{i\in S}x_i.
\end{equation}
\end{proposition}

The set functions $\alpha$ and $\beta$ that disjunctively and conjunctively define the polynomial function $p$ in Proposition~\ref{prop:lp dnf}
are not unique. For example, we have $x_1 \vee (x_1 \wedge x_2) = x_1 = x_1\wedge (x_1 \vee x_2)$. However, it can be shown \cite{Mar02c} that, from among all the possible set functions that disjunctively define a given lattice polynomial
function, only one is nondecreasing. Similarly, from among all the possible set functions that conjunctively define a given lattice polynomial
function, only one is nonincreasing. These particular set functions are given by $\alpha(S) = p(\mathbf{1}_S)$ and $\beta(S) =
p(\mathbf{1}_{[n]\setminus S})$ for all $S\subseteq [n]$, where $\mathbf{1}_S$ denotes the characteristic vector of $S$ in $\{0,1\}^n$. Thus, a
lattice polynomial function $p\colon\R^n\to\R$ can always be written as
$$
p(\bx)=\bigvee_{\textstyle{S\subseteq [n]\atop p(\mathbf{1}_S)=1}}\bigwedge_{i\in S}x_i=\bigwedge_{\textstyle{S\subseteq [n]\atop
p(\mathbf{1}_{[n]\setminus S})=0}}\bigvee_{i\in S}x_i.
$$

\begin{remark}
Now it becomes evident that any $n$-variable lattice polynomial function is a nondecreasing and continuous order invariant function in $\R^n$.
We will see in Proposition~\ref{prop:cindilp} that the converse is also true. A nondecreasing (or continuous) order invariant function in $\R^n$
is a lattice polynomial function.
\end{remark}

Denote by $p_{\alpha}^{\vee}$ (resp.\ $p_{\beta}^{\wedge}$) the lattice polynomial function disjunctively (resp.\ conjunctively) defined by a
given set function $\alpha$ (resp.\ $\beta$) as defined in Proposition~\ref{prop:lp dnf}. Let $f\colon\{0,1\}^n\to\{0,1\}$ be a nonconstant and
nondecreasing Boolean function. Then the lattice polynomial function $p_{\alpha}^{\vee}$, where $\alpha\colon 2^{[n]}\to\{0,1\}$ is defined by
$\alpha(S):=f(\mathbf{1}_S)$ for all $S\subseteq [n]$, is an extension to $\R^n$ of $f$. Indeed, we immediately have
$f(\mathbf{1}_S)=\alpha(S)=p_{\alpha}^{\vee}(\mathbf{1}_S)$ for all $S\subseteq [n]$. Consequently, any $n$-variable lattice polynomial function
is an extension to $\R^n$ of a nonconstant and nondecreasing Boolean function.

Throughout we will denote by $\cn$ the set of $\{0,1\}$-valued nonconstant and nondecreasing set functions on $[n]$. By definition, this set is
equipollent to the set of $n$-variable lattice polynomial functions, as well as to the set of nonconstant and nondecreasing Boolean
functions.\footnote{The problem of enumerating the number of distinct nondecreasing Boolean functions of $n$ variables is known as the
Dedekind's problem \cite{Kle69,KleMar75} (Sloane's integer sequence A000372). Although Dedekind first considered this question in 1897, there is
still no concise closed-form expression for this sequence.}

Now, regard the lattice polynomial function $p$ as a function from $E^n$ to $E$. If $E$ is a bounded lattice, we necessarily have $\bigvee_{x\in\varnothing}x:=\inf E$ and $\bigwedge_{x\in\varnothing}x:=\sup E$. Then, from (\ref{eq:pab}), we immediately see that $p\equiv\inf E$ if $\alpha\equiv 0$, and $p\equiv\sup E$ if $\alpha\equiv 1$. Thus we can
extend the definition of lattice polynomial functions by allowing the set function $\alpha$ to be constant.

Let $\cen$ denote the set $\cn$ completed with the constant set function $\alpha\equiv 0$, if $\inf E\in E$, and the constant set function
$\alpha\equiv 1$, if $\sup E\in E$. Evidently $\cen$ can be partially ordered by the standard partial order on set functions, namely $\alpha_1
\preccurlyeq \alpha_2$ if and only if $\alpha_1(S) \leqslant \alpha_2(S)$ for all $S\subseteq [n]$. We will refer to this partial order in the
subsequent sections.

\subsection{Special lattice polynomial functions}
\label{sec:plp}

We now consider the important special case of symmetric lattice polynomial functions. Denote by $x_{(1)},\ldots,x_{(n)}$ the {\em order
statistics}\/ resulting from reordering the variables $x_1,\ldots,x_n$ in nondecreasing order, that is, $x_{(1)}\leqslant\cdots\leqslant
x_{(n)}$. As Ovchinnikov~\cite[\S{7}]{Ovc96} observed, any order statistic is a symmetric lattice polynomial function. More precisely, for any
$k\in [n]$, we have
$$
x_{(k)}=\bigvee_{\textstyle{S\subseteq [n]\atop |S|=n-k+1}}\bigwedge_{i\in S}x_i=\bigwedge_{\textstyle{S\subseteq [n]\atop |S|=k}}\bigvee_{i\in
S}x_i.
$$
Conversely, Marichal~\cite[\S{2}]{Mar02c} showed that any symmetric lattice polynomial function is an order statistic.

Let us denote by $\mathrm{os}_k\colon \R^n\to\R$ the $k$th \emph{order statistic function}, that is, $\mathrm{os}_k(\bx):=x_{(k)}$. It is then easy to
see that, for any $S\subseteq [n]$, we have $\mathrm{os}_k(\mathbf{1}_S)=1$ if and only if $|S|\geqslant n-k+1$ and, likewise, we have
$\mathrm{os}_k(\mathbf{1}_{[n]\setminus S})=0$ if and only if $|S|\geqslant k$. Note that when $n$ is odd, $n = 2k-1$, the particular order statistic $x_{(k)}$ is the well-known \emph{median}\/ function
$$
{\rm median}(x_1,\ldots,x_{2k-1}) := x_{(k)}.
$$

Another special case of lattice polynomial functions is given by the projection functions, already introduced in \S\ref{sec:gfass3}. Recall
that, for any $k \in [n]$, the \emph{projection}\/ function ${\rm P}_k\colon\R^n \to \R$ associated with the $k$th argument is defined by ${\rm
P}_k(\bx) := x_k$. The projection function ${\rm P}_k$ consists in projecting $\bx \in \R^n$ onto the $k$th coordinate axis. As a particular aggregation function,
it corresponds to a dictatorial aggregation.

\subsection{Some aggregation properties}
\label{sec:fgrtre}

Lattice polynomial functions $p\colon E^n\to E$ are clearly continuous and nondecreasing functions. They are also order invariant functions in the
sense that they fulfill equation (\ref{eq:fppf3}) with arbitrary increasing bijections $\phi\colon E\to E$; see e.g.\ \cite{Ovc98c}.

Lattice polynomial functions also fulfill other properties shared by many aggregation functions. We now examine three of them : internality,
idempotency, and discretizability.

The most often encountered functions in aggregation theory are means or averaging functions, such as the weighted arithmetic means.
Cauchy~\cite{Cau21} considered in 1821 the mean of $n$ independent variables $x_1,\ldots,x_n$ as a function $F(x_1,\ldots,x_n)$ which should be
internal to the set of $x_i$ values.

\begin{definition}
$F\colon E^n \to \R$ is an \emph{internal} function if $\bigwedge_{i=1}^n x_i \leqslant F(\bx) \leqslant \bigvee_{i=1}^n x_i$ for all $\bx \in E^n$.
\end{definition}

Such means satisfy trivially the property of \emph{idempotency}, that is, if all $x_i$ are identical, $F(\bx)$ restitutes the common value.

\begin{definition}
$F\colon E^n \to \R$ is an \emph{idempotent} function if $F(x,\ldots,x) = x$ for all $x \in E$.
\end{definition}

Conversely, we can easily see that any nondecreasing and idempotent function $F\colon E^n \to \R$ is internal.

As any lattice polynomial function is clearly internal, it is a mean in the Cauchy sense. Thus, the internality property makes it possible to define
means even on ordinal scales (see, e.g.\ \cite{Ovc96}). For example, as a particular lattice polynomial function, the classical median
function (see \S\ref{sec:plp}), which gives the middle value of an odd-length sequence of ordered values, is a continuous, nondecreasing, and
symmetric mean defined on ordinal scales. To give a second example, consider the classical \emph{mode}\/ function, $\mathrm{mode}\colon E^n\to E$,
defined by\footnote{As usual, $\argmax$ stands for the argument of the maximum, that is to say, the value of the given argument for which the
value of the given expression attains its maximum value.}
\begin{equation}\label{eq:mode}
\mathrm{mode}(\bx):=\argmax_{r\in E}\,\sum_{i=1}^n \mathbf{1}_{\{0\}}(x_i-r),
\end{equation}
where the function $\mathbf{1}_{\{0\}}\colon\R\to\R$ is defined by $\mathbf{1}_{\{0\}}(0):=1$ and $\mathbf{1}_{\{0\}}(x):=0$ for all $x\neq 0$ (in
case of multiple values for $\argmax$, take the smallest one). This function, which gives the (lowest) most repeated value of a sequence of
values, is a symmetric mean defined on ordinal scales, and even on nominal scales.\footnote{The admissible transformations associated with a
nominal scale are one-to-one transformations (injections) of $E$ into itself; see \cite[p.~66]{Rob79}.} However, since the mode function
is not nondecreasing, it is not a lattice polynomial function.

We can also observe that any lattice polynomial function is \emph{discretizable}\/ in the sense that it always yields the value of one of its variables.
This property was actually introduced in the framework of triangular norms (see e.g.\ \cite{DeBMes03,Fod00}) but is easily extended to any function as follows.

\begin{definition}\label{de:disc}
$F\colon E^n \to E$ is a \emph{discretizable} function if $F(\bx)\in\{x_1,\ldots,x_n\} \cup B[E]$ for all $\bx\in E^n$.
\end{definition}

We can readily prove~\cite{Mes01} that $F\colon E^n \to E$ is a discretizable function if and only if, for any nonempty finite subset $C\subset E$
and any $\bx\in (C\cup B[E])^n$, we have $F(\bx)\in C\cup B[E]$. Thus, this property means that the domain and range of $F$ can be restricted to
a finite or countable chain.

Another interesting property is {\em self-duality} (for bounded $E$, see \cite{GarMar08} and the
references therein), which is fulfilled for example by the median and the mode functions.

\begin{definition}\label{de:selfduality}
Let $\psi\colon E\to E$ be a decreasing and involutive (i.e., $\psi\circ\psi=\mathrm{id}$) bijection (hence necessarily $B[E]$ is not a singleton).
\begin{itemize}
\item The {\em $\psi$-dual} of a function $F\colon E^n\to E$ is the function $F_{\psi}\colon E^n\to E$, defined by
$$
F_{\psi}(\bx):=\psi^{-1}\big(F\big(\psi(x_1),\ldots,\psi(x_n)\big)\big).
$$
\item A function $F\colon E^n\to E$ is said to be {\em $\psi$-self-dual} if $F_{\psi}=F$.

\item A function $F\colon E^n\to E$ is said to be {\em weakly self-dual} if it is $\psi$-self-dual for some decreasing and involutive bijection
$\psi\colon E\to E$.
\end{itemize}
If $E$ is bounded, then the only affine decreasing bijection from $E$ onto itself is given by $\psi^d(x):=\inf E+\sup E-x$, and $\psi^d$-duality
is then called \emph{duality}, with notation $F^d:=F_{\psi^d}$. A function $F\colon E^n\to E$ is said to be {\em self-dual} if $F^d=F$.
\end{definition}

\begin{remark}
\begin{enumerate}
\item[(i)] By definition, for any function $F\colon E^n\to E$ and any decreasing and involutive bijection $\psi\colon E\to E$, we have
$(F_{\psi})_{\psi}=F$.

\item[(ii)] We note that $\psi$-duality is an example of {\em $\psi$-conjugacy} \cite[Chapter~8]{KucChoGer90}, whose definition is the same
except that it does not require $\psi$ to be decreasing nor involutive. We also observe that the classic notion of duality in ordered sets
concerns the converse order relations; see \cite[p.~3]{Bir67}.
\end{enumerate}
\end{remark}

Assume that $B[E]$ is not a singleton and let $\alpha\in\cn$. We can straightforwardly show that the lattice polynomial function
$p_{\alpha}^{\vee}\colon E^n\to E$ is weakly self-dual if and only if $\alpha^d=\alpha$, where $\alpha^d\in\cn$ is the dual of $\alpha$, defined by
$\alpha^d(S):=1-\alpha([n]\setminus S)$.

The special case of order statistics is dealt with in the next immediate result (see \cite[\S{5}]{Mar02c}), which characterizes the
median as the only weakly self-dual order statistic.

\begin{proposition}\label{prop:charmed}
Assume that $n$ is odd and that $B[E]$ is not a singleton. The $k$th order statistic function ${\rm os}_k\colon E^n\to E$ is weakly self-dual if and
only if $n=2k-1$. In this case ${\rm os}_k$ is the median function.
\end{proposition}

%The following description of $\psi$-self-dual lattice polynomial functions is also worth mentioning (for bounded $E$, see Yanovskaya
%\cite[p.~827]{Yan89}): (to check!!!)
%
%\begin{proposition}
%Assume $B[E]$ is not a singleton and let $\psi:E\to E$ be a decreasing bijection. A lattice polynomial $p:E^n\to E$ is $\psi$-self-dual if and
%only if there exist lattice polynomials $p':E^n\to E$ and $p'':E^n\to E$, with $p''$ $\psi$-self-dual, such that
%$$
%p(\bx)={\rm median}(p'(\bx),p'_{\psi}(\bx),p''(\bx)).
%$$
%\end{proposition}

%---------------------------------------------------------------------------------------------- Section 4
\section{Order invariant functions}

The first meaningful aggregation functions we consider are the \emph{order invariant functions}, which were first investigated (as ordinally
stable functions) by Marichal and Roubens~\cite{MarRou93}, and then by many other authors; see
\cite{BarDre04,FodRou95,KamOvc95,Mar98,Mar02c,MarMat01,MarMes04,MarMesRuc05,Mes01,MesRuc04,Ovc96,Ovc98c,OvcDuk00,OvcDuk02}.

\subsection{Definition and first results}

Let $x_1,\ldots,x_n$ be independent variables defining the same ordinal scale, with domain $E$, and suppose that, when aggregating these
variables by a function $F\colon E^n\to E$, we require that the dependent variable
\begin{equation}\label{eq:yem}
x_{n+1}=F(x_1,\ldots,x_n)
\end{equation}
defines the same scale. As equation (\ref{eq:yem}) should represent a meaningful relation between the independent and dependent variables, the
aggregation function $F$ should be invariant under actions from $\Phi[E]$. That is, $\phi(x_{n+1})=F(\phi(x_1),\ldots,\phi(x_n))$ for all $\phi\in \Phi[E]$. Thus, the order invariance property is defined as follows.

\begin{definition}\label{de:inv}
$F\colon E^n \to E$ is said to be an \emph{order invariant}\/ function if
$$
F\big(\phi(x_1),\ldots,\phi(x_n)\big) = \phi\big(F(x_1,\ldots,x_n)\big)
$$
for all $\bx\in E^n$ and all $\phi \in \Phi[E]$.
\end{definition}

The following result (see Propositions~\ref{prop:eondin} and \ref{prop:eocoin} below) shows that the lattice polynomial functions are the most
prominent order invariant functions (however see Theorem~\ref{thm:nvkjfd} for a full description of order invariant functions).

\begin{proposition}\label{prop:cindilp}
Assume that $E$ is open and consider a function $F\colon E^n\to E$. Then the following three assertions are equivalent:
\begin{enumerate}
\item[(i)] $F$ is a nondecreasing order invariant function.

\item[(ii)] $F$ is a continuous order invariant function.

\item[(iii)] $F$ is a lattice polynomial function.
\end{enumerate}
\end{proposition}

Proposition~\ref{prop:cindilp} poses the interesting question of how we can interpret the continuity property for order invariant functions. Let
$\Phi'[E]$ be the superset of $\Phi[E]$ consisting of the continuous nondecreasing surjections $\phi\colon E\to E$. The following result
\cite[\S{5.2}]{MarMes04}, inspired from \cite[Proposition 2]{BouPir97}, shows that the conjunction of continuity and order invariance is
equivalent to requiring that the admissible transformations belong to $\Phi'[E]$.

\begin{proposition}%\label{prop:cif/ap}
$F\colon E^n\to E$ is a continuous order invariant function if and only if $F(\phi(x_1),\ldots,\phi(x_n)) = \phi(F(x_1,\ldots,x_n))$ for all $\bx\in E^n$ and all $\phi \in \Phi'[E]$.
\end{proposition}

Let $\Phi''[E]$ be the superset of $\Phi[E]$ consisting of all the monotone bijections of $E$ onto itself (assuming that $B[E]$ is not a
singleton). It is clear (for bounded $E$, see e.g.\ \cite[\S{3}]{MesRuc04}) that the conjunction of weak
self-duality (cf.\ Definition~\ref{de:selfduality}) and order invariance is equivalent to requiring that the admissible transformations belong
to $\Phi''[E]$. The independent and dependent variables then define a \emph{nominal} scale.

\begin{proposition}%\label{prop:cif/ap}
Assume that $B[E]$ is not a singleton. $F\colon E^n\to E$ is a  weakly self-dual order invariant function if and only if $F(\phi(x_1),\ldots,\phi(x_n)) = \phi(F(x_1,\ldots,x_n))$ for all $\bx\in E^n$ and all $\phi \in \Phi''[E]$.
\end{proposition}

\subsection{General descriptions}

When $E$ is open we have the following description (see \cite[Theorem~5.1]{Ovc98c}).

\begin{proposition}\label{prop:eoin}
Assume that $E$ is open. Then $F\colon E^n\to E$ is an order invariant function if and only if there exists a mapping $\xi\colon\ien\to [n]$ such that
$F|_I={\rm P}_{\xi(I)}|_I$ for all $I\in\ien$.
\end{proposition}

This result shows that, when $E$ is open, the restriction of $F$ to any minimal order invariant subset is a projection function onto one
coordinate. That is, for any $I\in\ien$, there exists $k_I\in [n]$ such that $F|_I=\mathrm{P}_{k_I}|_I$. Clearly, such a function is internal
and hence idempotent.

As an example, any nonconstant lattice polynomial function is a continuous, nondecreasing, idempotent, and order invariant function. On the
other hand, the mode function (\ref{eq:mode}) is an idempotent and order invariant function that is neither continuous nor nondecreasing.

When $E$ is not open, the restriction of $F$ to any minimal order invariant subset reduces to a constant function or a projection function onto
one coordinate (see \cite{MarMesRuc05,MesRuc04}).

\begin{theorem}\label{thm:nvkjfd}
$F\colon E^n\to E$ is an order invariant function if and only if there exists a mapping $\xi\colon\ien\to [n]$ such that, for any $I\in\ien$,
\begin{itemize}
\item either $F|_I\equiv c\in B[E]$ (assuming $B[E]\neq\varnothing$),

\item or $F|_I={\rm P}_{\xi(I)}|_I$.
\end{itemize}
\end{theorem}

\begin{remark}
It was proved in \cite[Proposition~3.1]{Mar02c} (see \cite{KamOvc95,Ovc96,Ovc98c} for preliminary results), that any order invariant function is
discretizable, and hence it is internal whenever $E$ is open; it is clear that this follows from Theorem~\ref{thm:nvkjfd}. For instance, the
mode function (\ref{eq:mode}) is order invariant and hence discretizable. The converse is not true. For example, the function $F\colon (0,1)^2\to
(0,1)$ defined by
$$
F(x_1,x_2):=
\begin{cases}
x_1, & \mbox{if $x_1+x_2<1$},\\
x_2, & \mbox{otherwise},
\end{cases}
$$
is discretizable but not order invariant.
\end{remark}

When an order invariant function is idempotent, clearly it must be a projection function on the open diagonal of $E^n$, also on $I=\{(\inf
E,\ldots,\inf E)\}$ (if $\inf E\in E$), and on $I=\{(\sup E,\ldots,\sup E)\}$ (if $\sup E\in E$).

\subsection{The nondecreasing case}

We now present descriptions of order invariant functions which are nondecreasing. The following result (see \cite[Corollary~4.4]{Mar02c}) shows that, when $E$ is open, the family of nondecreasing order invariant functions in $E^n$ is identical
to that of lattice polynomial functions in $E^n$.

\begin{proposition}\label{prop:eondin}
Assume that $E$ is open. Then $F\colon E^n\to E$ is a nondecreasing order invariant function if and only if it is a lattice polynomial function.
\end{proposition}

\begin{corollary}\label{cor:gdfg}
Assume that $E$ is open. Then $F\colon E^n\to E$ is a symmetric, nondecreasing, and order invariant function if and only if it is an order statistic
function.
\end{corollary}

Combining Proposition~\ref{prop:charmed} with Corollary~\ref{cor:gdfg} immediately yields the following axiomatization of the median function.

\begin{corollary}
Assume that $n$ is odd and that $E$ is open. Then $F\colon E^n\to E$ is a symmetric, weakly self-dual, nondecreasing, and order invariant function if
and only if it is the median function.
\end{corollary}

A complete description of nondecreasing order invariant functions in $E^n$, with open or non-open interval $E$, is given in the following
theorem (see \cite{MarMesRuc05,MesRuc04}). It shows that discontinuities of $F$ may
occur only on the border of $E^n$. Recall that the lattice polynomial function in $E^n$ disjunctively defined by $\alpha\in\cen$ is denoted
$p_{\alpha}^{\vee}$ (see \S\ref{sec:sfdaag}).

\begin{theorem}
$F\colon E^n\to E$ is a nondecreasing order invariant function if and only if there exists a nondecreasing mapping $\xi\colon\isen\to\cen$ such that
$F|_I=p_{\xi(I)}^{\vee}|_I$ for all $I\in\isen$.
\end{theorem}

\begin{example}
Consider the semiopen interval $E:=[a,b)$.\footnote{Here the poset $\isen$ contains 8 elements (a point, three open line segments, three open
square facets, and an open cube).} The function $F\colon [a,b)^3\to [a,b)$ defined by
$$
F(x_1,x_2,x_3):=\begin{cases} a, & \mbox{if $x_1=a$}, \\ x_3, & \mbox{if $x_1\neq a$ and $x_2=a$}, \\ x_1\vee x_2\vee x_3, &
\mbox{otherwise},\end{cases}
$$
is a nondecreasing order invariant function in $[a,b)^3$.
\end{example}

\begin{corollary}
$F\colon E^n\to E$ is a nondecreasing, idempotent, and order invariant function if and only if there exists a nondecreasing mapping
$\xi\colon\isen\to\cen$, where $\xi[(E^\circ)^n]$ is nonconstant, such that $F|_I=p_{\xi(I)}^{\vee}|_I$ for all $I\in\isen$.
\end{corollary}

\subsection{The continuous case}\label{sec:OIcont}

We now consider the family of continuous order invariant functions. It was shown in \cite[Corollary~4.2]{Mar02c} that, when $E$ is open, this
family is identical to the family of lattice polynomial functions in $E^n$; see also \cite[\S{3.4.2}]{Mar98}.

\begin{proposition}\label{prop:eocoin}
Assume that $E$ is open. Then $F\colon E^n\to E$ is a continuous order invariant function if and only if it is a lattice polynomial function.
\end{proposition}

\begin{remark}
Note that this result was independently stated and proved earlier by Ovchinnikov~\cite[Theorem~5.3]{Ovc98} in the more general setting where the
range of variables is a doubly homogeneous simple order, that is, a simple order $X$ satisfying the following property:
\begin{quote}
For any $x_1, x_2, y_1, y_2 \in X$, with $x_1 < x_2$ and $y_1 < y_2$, there is an automorphism $\phi\colon X \to X$ such that $\phi(x_1) = y_1$ and
$\phi(x_2) = y_2$.
\end{quote}
As any open interval $E$ of the real line is clearly a doubly homogeneous simple order, Ovchinnikov's result encompasses that of
Proposition~\ref{prop:eocoin}.\footnote{Note that the extension of this result to the (infinite) case of functional operators was described in \cite{OvcDuk00}; see also~\cite{OvcDuk02}.}
\end{remark}

A complete description of continuous order invariant function in $E^n$ was stated in \cite[Corollary~4.3]{Mar02c} as follows (see
also \cite{MarMesRuc05}).

\begin{theorem}\label{thm:co/inv}
$F\colon E^n\to E$ is a continuous order invariant function if and only if there exists $\alpha\in\cen$ such that $F=p_{\alpha}^{\vee}$.
\end{theorem}

Theorem \ref{thm:co/inv} actually says that a continuous order invariant function $F\colon E^n\to E$ is either the constant function $F\equiv\inf E$
if $\inf E\in E$, or the constant function $F\equiv\sup E$ if $\sup E\in E$, or any lattice polynomial function in $E^n$ (any order statistic
function in $E^n$ if $F$ is symmetric).

\begin{remark}
From Theorems~\ref{thm:nvkjfd} and \ref{thm:co/inv} it follows that a function $F\colon E^n\to E$ is a lattice polynomial function if and only if its
restriction to each closed simplex of the standard triangulation of $E^n$ is a projection function onto one coordinate (see also \cite[Proposition~2.1]{Mar02c}).
\end{remark}

\begin{corollary}\label{cor:coidin}
$F\colon E^n\to E$ is a continuous, idempotent, and order invariant function if and only if it is a lattice polynomial function.
\end{corollary}

\begin{corollary}\label{cor:sycoidin}
$F\colon E^n\to E$ is a symmetric, continuous, idempotent, and order invariant function if and only if it is an order statistic function.
\end{corollary}

\begin{corollary}\label{cor:dsfs}
Assume that $n$ is odd and that $B[E]$ is not a singleton. Then $F\colon E^n\to E$ is a symmetric, weakly self-dual, continuous, idempotent, and order
invariant function if and only if it is the median function.
\end{corollary}

\begin{remark}
By combining Proposition~\ref{prop:charmed} with Theorem \ref{thm:co/inv}, we immediately see that idempotency is not necessary in
Corollary~\ref{cor:dsfs}. Indeed, a weak self-dual lattice polynomial function cannot be constant and hence it is necessarily idempotent.
\end{remark}

%---------------------------------------------------------------------------------------------- Section 5
\section{Comparison meaningful functions on a single ordinal scale}

We now present the class of \emph{comparison meaningful functions on a single ordinal scale}. These functions were introduced first by
Orlov~\cite{Orl81} and then investigated by many other authors; see
\cite{KamOvc95,Mar98,Mar01,Mar02c,MarMat01,MarMes04,MarMesRuc05,Ovc96,Yan89}.

\subsection{Definition and first results}

Let $x_1,\ldots,x_n$ be independent variables defining the same ordinal scale, with domain $E$, and suppose that, when aggregating these
variables by a function $F\colon E^n\to \R$, we require that the dependent variable $x_{n+1}=F(x_1,\ldots,x_n)$ defines an ordinal scale, with an
arbitrary domain in $\R$. According to Luce's principle \cite{Luc59}, any admissible transformation of the independent variables must lead to an
admissible transformation of the dependent variable. This condition can be formulated as follows.

\begin{definition}\label{de:cmfsos}
$F\colon E^n \to \R$ is said to be a \emph{comparison meaningful function on a single ordinal scale}\/ if, for any $\bphi\in \Phi_n[E]$, there is a
strictly increasing mapping $\psi_{\bphi}\colon \ran(F)\to\ran(F)$ such that $F[\bphi(\bx)]=\psi_{\bphi}[F(\bx)]$ for all $\bx\in E^n$.
\end{definition}

Comparison meaningful functions on a single ordinal scale were first introduced by Orlov~\cite{Orl81} as those functions preserving the
comparison of aggregated values when changing the scale defined by the independent variables.\footnote{A general study on meaningfulness of
ordinal comparisons can be found in \cite{RobRos94}.} We paraphrase from Orlov:
\begin{quote}
When one compares two sets of objects according to a criterion, it is sometimes required to evaluate each object on the same ordinal scale
(e.g., by means of measurement or expert estimate). The aggregated values of the evaluations corresponding to each set of objects are computed
by a certain aggregation function, and then compared together. It is natural to require that the inferences made from this comparison are
meaningful, that is, depend only on the initial information, but not on the scale used.\footnote{More generally, a statement using scales of
measurement is said to be {\sl meaningful}\/ if its truth or falsity is invariant when every scale is replaced by another acceptable version of
it; see \cite[p.~59]{Rob79}.}
\end{quote}
The equivalence between Definition~\ref{de:cmfsos} and Orlov's definition can be formulated mathematically as follows (see \cite{MarMes04}).

\begin{proposition}\label{prop:cmim}
$F\colon E^n \to \R$ is a comparison meaningful function on a single ordinal scale if and only if
$$
F(\bx)\loe F(\bx') \quad\Rightarrow\quad F\big(\bphi(\bx)\big)\loe F\big(\bphi(\bx')\big)
$$
for all $\bx,\bx'\in E^n$ and all $\bphi\in\Phi_n[E]$.\footnote{Equivalently, $F\colon E^n \to \R$ is a comparison meaningful function on a single
ordinal scale if and only if
$$
F(\bx)\leqslant F(\bx') \quad\Leftrightarrow\quad F\big(\bphi(\bx)\big)\leqslant F\big(\bphi(\bx')\big)
$$
for all $\bx,\bx'\in E^n$ and all $\bphi\in\Phi_n[E]$.}
\end{proposition}

\begin{remark}
Although the condition in Proposition~\ref{prop:cmim} is natural and even mandatory to aggregate ordinal values, it severely restricts the
allowable operations for defining a meaningful aggregation function. For example, the comparison of two arithmetic means is meaningless on an
ordinal scale. Indeed, considering the pairs of values $(3,5)$ and $(1,8)$, we have $\frac 12(3+5) <\frac 12(1+8)$ and, using any admissible
transformation $\phi$ such that $\phi(1)=1$, $\phi(3)=4$, $\phi(5)=7$, and $\phi(8)=8$, we have $\frac 12(\phi(3)+\phi(5))
>\frac 12(\phi(1)+\phi(8))$.
\end{remark}

Order invariant functions and comparison meaningful functions on a single ordinal scale can actually be related through the idempotency
property. Indeed, when a comparison meaningful function on a single ordinal scale is idempotent then the output scale must coincide with the
input scale. This result is stated in the next proposition (see \cite[Proposition~3.3]{Mar02c} and preliminary work in \cite{KamOvc95,Ovc96}).

\begin{proposition}\label{prop:idcmin}
Consider a function $F\colon E^n \to E$.
\begin{itemize}
\item If $F$ is idempotent and comparison meaningful on a single ordinal scale then it is order invariant.

\item If $F$ is order invariant then it is comparison meaningful on a single ordinal scale.

\item If $E$ is open then $F$ is idempotent and comparison meaningful on a single ordinal scale if and only if it is order invariant.
\end{itemize}
\end{proposition}

Just as for order invariant functions, continuity of comparison meaningful functions on a single ordinal scale can be interpreted by means of
the set $\Phi'[E]$ of continuous nondecreasing surjections from $E$ onto $E$; see \cite[\S{5.2}]{MarMes04}. Denote by
$\Phi'_n[E]$ the diagonal restriction of $\Phi'[E]^n$ (see \S\ref{sec:gfass3}).

\begin{proposition}
$F\colon E^n \to \R$ is a continuous and comparison meaningful function on a single ordinal scale if and only if, for any $\bphi\in \Phi'_n[E]$,
there is a continuous and nondecreasing mapping $\psi_{\phi}\colon\ran(F)\to\ran(F)$ such that $F\big(\bphi(\bx)\big)=\psi_{\bphi}\big(F(\bx)\big)$
for all $\bx\in E^n$.
\end{proposition}

\subsection{General descriptions}

The class of comparison meaningful functions on a single ordinal scale can be described as follows (see \cite[Theorem~3.1]{MarMesRuc05}).

\begin{theorem}
$F\colon E^n\to\R$ is a comparison meaningful function on a single ordinal scale if and only if, for any $I\in\ien$, there exist an index $k_I\in [n]$
and a strictly monotonic or constant function $g_I\colon \mathrm{P}_{k_I}(I)\to\R$ such that $F|_I=(g_I\circ\mathrm{P}_{k_I})|_I$, where, for any
$I,I'\in\ien$,
\begin{itemize}
\item either $g_I=g_{I'}$,

\item or $\ran(g_I)=\ran(g_{I'})$ is a singleton,

\item or $\ran(g_I)<\ran(g_{I'})$,

\item or $\ran(g_I)>\ran(g_{I'})$.\footnote{Note that $\ran(g_I)<\ran(g_{I'})$ means that for all $r\in\ran(g_I)$ and all $r'\in\ran(g_{I'})$,
we have $r<r'$.}
\end{itemize}
\end{theorem}

Thus, a comparison meaningful function on a single ordinal scale reduces, on each minimal order invariant subset of $E^n$, to a constant or a
transformed projection function onto one coordinate.

\begin{example}
We have seen in Example~\ref{ex:fggger} that there are eleven minimal order invariant subsets in the unit square $[0,1]^2$, namely
\begin{itemize}
\item $I_1:=\{(0,0)\}$, $I_2:=\{(1,0)\}$, $I_3:=\{(1,1)\}$, $I_4:=\{(0,1)\}$,

\item $I_5:=(0,1)\times\{0\}$, $I_6:=\{1\}\times (0,1)$, $I_7:=(0,1)\times\{1\}$, $I_8:=\{0\}\times (0,1)$,

\item $I_9:=\{(x_1,x_2)\mid 0<x_1=x_2 <1\}$, $I_{10}:=\{(x_1,x_2)\mid 0<x_1<x_2<1\}$, $I_{11}:=\{(x_1,x_2)\mid 0<x_2<x_1<1\}$.
\end{itemize}
Let $k_{I_j}:=1$ and $g_{I_j}(x):=1-x$ if $j\in\{1,2,3,5,6,9,11\}$, and $k_{I_j}:=2$ and $g_{I_j}(x):=2x-3$ if $j\in\{4,7,8,10\}$, where always
$x\in {\rm P}_{k_{I_j}}(I_j)$. Then the corresponding comparison meaningful function $F\colon [0,1]^2\to \R$ is given by
$$
F(x_1,x_2):= \begin{cases} 1-x_1, & \mbox{if $x_1\geqslant x_2$}, \\ 2x_2-3, & \mbox{if $x_1<x_2$}.\end{cases}
$$
\end{example}

When a comparison meaningful function on a single ordinal scale is idempotent, it must satisfy $g_I(x)=x$, for all $x\in {\rm P}_{k_I}(I)$,
whenever either $I$ is the open diagonal of $E^n$, or $I=\{(\inf E,\ldots,\inf E)\}$ (if $\inf E\in E$), or $I=\{(\sup E,\ldots,\sup E)\}$ (if
$\sup E\in E$).

\subsection{The nondecreasing case}

The following result~\cite{MarMesRuc05} yields, when $E$ is open, a description of all nondecreasing comparison meaningful functions $F\colon E^n\to
\R$ on a single ordinal scale.

\begin{proposition}\label{prop:eondcm}
Assume that $E$ is open. Then $F\colon E^n\to\R$ is a nondecreasing comparison meaningful function on a single ordinal scale if and only if there
exist $\alpha\in\cn$ and a strictly increasing or constant function $g\colon E\to\R$ such that $F=g\circ p_{\alpha}^{\vee}$.
\end{proposition}

As we can see, all the functions described in Proposition~\ref{prop:eondcm} are continuous up to possible discontinuities of the function $g$.

The following corollaries~\cite[Theorem~4.4]{Mar02c} (see~\cite[Theorem~3.1]{MarMat01} for preliminary results) immediately follow from
Proposition~\ref{prop:eondcm}.

\begin{corollary}
Assume that $E$ is open. Then $F\colon E^n\to \R$ is a nondecreasing, idempotent, and comparison meaningful function on a single ordinal scale if and
only if it is a lattice polynomial function.
\end{corollary}

\begin{corollary}\label{cor:sdsds}
Assume that $E$ is open. Then $F\colon E^n\to \R$ is a symmetric, nondecreasing, idempotent, and comparison meaningful function on a single ordinal
scale if and only if it is an order statistic function.
\end{corollary}

\begin{corollary}
Assume that $n$ is odd and that $E$ is open. Then $F\colon E^n\to \R$ is a symmetric, weakly self-dual, nondecreasing, idempotent, and comparison
meaningful function on a single ordinal scale if and only if it is the median function.
\end{corollary}

A complete description of nondecreasing comparison meaningful functions $F\colon E^n\to \R$ on a single ordinal scale is given in the next
theorem~\cite[Corollary~4.1]{MarMesRuc05}. Let $\mathcal{G}(E)$ be the set of all strictly increasing or constant real functions $g$ defined either
on $E^\circ$, or on the singleton $\{\inf E\}\cap E$, or on the singleton $\{\sup E\}\cap E$ (if these singletons exist). This set is partially
ordered as follows: $g_1 \preccurlyeq g_2$ if either $g_1=g_2$, or $\ran(g_1)=\ran(g_2)$ is a singleton, or $\ran(g_1)<\ran(g_2)$.

\begin{theorem}
$F\colon E^n\to \R$ is a nondecreasing comparison meaningful function on a single ordinal scale if and only if there exist nondecreasing mappings
$\gamma\colon \isen\to \mathcal{G}(E)$ and $\xi\colon \isen\to\cen$ such that $F|_I=(\gamma(I)\circ p_{\xi(I)}^{\vee})|_I$ for all $I\in\isen$.
\end{theorem}

If furthermore $F$ is idempotent, then by nondecreasing monotonicity, we have $\ran(F)=E$ and, by Proposition~\ref{prop:idcmin}, $F$ is order
invariant. Hence we have the following corollary.

\begin{corollary}
$F\colon E^n\to \R$ is a nondecreasing, idempotent, and comparison meaningful function on a single ordinal scale if and only if there exists a
nondecreasing mapping $\xi\colon \isen\to\cen$, where $\xi[(E^\circ)^n]$ is nonconstant, such that $F|_I=p_{\xi(I)}^{\vee}|_I$ for all $I\in\isen$.
\end{corollary}

\subsection{The continuous case}

Based on a preliminary result~\cite[\S{4}]{MarMat01} (see also~\cite[\S{3.4.2}]{Mar98}), a full description of continuous comparison meaningful
functions on a single ordinal scale was given in \cite[Theorem~4.2]{Mar02c} as follows.

\begin{theorem}
$F\colon E^n\to \R$ is a continuous comparison meaningful function on a single ordinal scale if and only if there exist $\alpha\in\cn$ and a
continuous and strictly monotonic or constant function $g\colon E\to\R$ such that $F=g\circ p_{\alpha}^{\vee}$.
\end{theorem}

\begin{corollary}\label{cor:coidcm}
$F\colon E^n\to \R$ is a continuous, idempotent, and comparison meaningful function on a single ordinal scale if and only if it is a lattice
polynomial function.
\end{corollary}

\begin{remark}
The result in Corollary~\ref{cor:coidcm} was stated and proved first in social choice theory by Yanovskaya~\cite[Theorem~1]{Yan89} when $E=\R$.
\end{remark}

\begin{corollary}
$F\colon E^n\to \R$ is a symmetric, continuous, and comparison meaningful function on a single ordinal scale if and only if there exist $k\in [n]$ and
a continuous strictly monotonic or constant function $g\colon E\to\R$ such that $F=g\circ {\rm os}_k$.
\end{corollary}

\begin{corollary}\label{cor:sycoidcm}
$F\colon E^n\to \R$ is a symmetric, continuous, idempotent, and comparison meaningful function on a single ordinal scale if and only if it is an order
statistic function.
\end{corollary}

\begin{remark}
A slightly stronger version of the result in Corollary~\ref{cor:sycoidcm}, consisting in replacing idempotency with internality, was actually
proved first by Orlov~\cite{Orl81} in $\R^n$, then by Marichal and Roubens~\cite[Theorem~1]{MarRou93} in $E^n$ (see also
\cite[Theorem~3.4.13]{Mar98}), and finally by Ovchinnikov~\cite[Theorem~4.3]{Ovc96} in the more general framework where the range of variables
is a simple order $X$ whose open intervals are homogeneous and nonempty (see also~\cite[\S{6}]{Ovc98}).
\end{remark}

\begin{corollary}
Assume that $n$ is odd and that $B[E]$ is not a singleton. Then $F\colon E^n\to \R$ is a symmetric, weakly self-dual, continuous, idempotent, and
comparison meaningful function on a single ordinal scale if and only if it is the median function.
\end{corollary}

%---------------------------------------------------------------------------------------------- Section 6
\section{Comparison meaningful functions on independent ordinal scales}

In this section we present the class of \emph{comparison meaningful functions on independent ordinal scales}, which were introduced by Acz\'el
and Roberts~\cite[Case~\#21]{AczRob89} and studied by Kim~\cite{Kim90} (see preliminary work in Osborne~\cite{Osb70}) and then investigated by
some other authors; see \cite{Mar02c,MarMat01,MarMes04,MarMesRuc05}.

\subsection{Definition and first results}\label{sec:sdsgf}

Let $x_1,\ldots,x_n$ be independent variables defining independent ordinal scales, with a common domain $E$, and suppose that, when aggregating
these variables by a function $F\colon E^n\to \R$, we require that the dependent variable $x_{n+1}=F(x_1,\ldots,x_n)$ defines an ordinal scale, with
an arbitrary domain in $\R$. This condition can be formulated as follows.

\begin{definition}\label{de:cmfios}
$F\colon E^n \to \R$ is said to be a \emph{comparison meaningful function on independent ordinal scales}\/ if, for any $\bphi\in \Phi[E]^n$, there
is a strictly increasing mapping $\psi_{\bphi}\colon \ran(F)\to\ran(F)$ such that $F\big(\bphi(\bx)\big)=\psi_{\bphi}\big(F(\bx)\big)$ for all $\bx\in
E^n$.
\end{definition}

Comparison meaningful functions on independent ordinal scales can also be defined as those functions preserving the comparison of aggregated
values when changing the scales defined by the independent variables (see \cite{MarMes04}).

\begin{proposition}\label{prop:cmim2}
$F\colon E^n \to \R$ is a comparison meaningful function on independent ordinal scales if and only if
$$
F(\bx)\loe F(\bx') \quad\Rightarrow\quad F\big(\bphi(\bx)\big)\loe F\big(\bphi(\bx')\big)
$$
for all $\bx,\bx'\in E^n$ and all $\bphi\in\Phi[E]^n$.\footnote{Equivalently, $F\colon E^n \to \R$ is a comparison meaningful function on
independent ordinal scales if and only if
$$
F(\bx)\leqslant F(\bx') \quad\Leftrightarrow\quad F\big(\bphi(\bx)\big)\leqslant F\big(\bphi(\bx')\big)
$$
for all $\bx,\bx'\in E^n$ and all $\bphi\in\Phi[E]^n$.}
\end{proposition}

Comparison meaningfulness on independent ordinal scales is a very strong condition, much stronger than comparison meaningfulness on a single
ordinal scale. For example, it was proved~\cite[Lemma~5.2]{Mar02c} that this condition reduces any lattice polynomial function to a projection
function onto one coordinate.

Regarding continuity of comparison meaningful function on independent ordinal scales, it can be interpreted in the same way as for comparison
meaningful function on a single ordinal scale; see \cite[\S{5.2}]{MarMes04}. Consider again the set $\Phi'[E]$ of continuous
nondecreasing surjections from $E$ onto $E$.

\begin{proposition}
$F\colon E^n \to \R$ is a continuous and comparison meaningful function on independent ordinal scales if and only if, for any $\bphi\in \Phi'[E]^n$,
there is a continuous and nondecreasing mapping $\psi_{\bphi}\colon \ran(F)\to\ran(F)$ such that $F\big(\bphi(\bx)\big)=\psi_{\bphi}\big(F(\bx)\big)$
for all $\bx\in E^n$.
\end{proposition}

\subsection{General descriptions}

The description of comparison meaningful functions on independent ordinal scales is very similar to that of comparison meaningful functions on a
single ordinal scale. The result can be formulated as follows~\cite[Corollary~3.1]{MarMesRuc05}.

\begin{theorem}\label{thm:fdgfasd}
$F\colon E^n\to\R$ is a comparison meaningful function on independent ordinal scales if and only if, for any $I\in\isen$, there exist an index $k_I\in
[n]$ and a strictly monotonic or constant function $g_I\colon \mathrm{P}_{k_I}(I)\to\R$ such that $F|_I=(g_I\circ\mathrm{P}_{k_I})|_I$, where, for any
$I,I'\in\isen$,
\begin{itemize}
\item either $g_I=g_{I'}$,

\item or $\ran(g_I)=\ran(g_{I'})$ is a singleton,

\item or $\ran(g_I)<\ran(g_{I'})$,

\item or $\ran(g_I)>\ran(g_{I'})$.
\end{itemize}
\end{theorem}

Thus, a comparison meaningful function on independent ordinal scales reduces, on each minimal strongly order invariant subset of $E^n$, to a
constant or a transformed projection function onto one coordinate.

When a comparison meaningful function on independent ordinal scales is idempotent, it must satisfy $g_I(x)=x$, for all $x\in {\rm P}_{k_I}(I)$,
whenever either $I=(E^\circ)^n$, or $I=\{(\inf E,\ldots,\inf E)\}$ (if $\inf E\in E$), or $I=\{(\sup E,\ldots,\sup E)\}$ (if $\sup E\in E$).

When $E$ is open, the family $\isen$ reduces to $\{E^\circ\}$, thus considerably simplifying Theorem~\ref{thm:fdgfasd} as follows.

\begin{proposition}\label{prop:eoscm}
Assume that $E$ is open. Then $F\colon E^n\to \R$ is a comparison meaningful function on independent ordinal scales if and only if there exist $k\in
[n]$ and a strictly monotonic or constant function $g\colon E\to\R$ such that $F=g\circ {\rm P}_k$.
\end{proposition}

\begin{corollary}
Assume that $E$ is open. Then $F\colon E^n\to \R$ is an idempotent and comparison meaningful function on independent ordinal scales if and only if it
is a projection function.
\end{corollary}

It follows from Proposition~\ref{prop:eoscm} that, when $E$ is open and $n\geqslant 2$, any symmetric and comparison meaningful function on
independent ordinal scales is necessarily a constant function. In this case, it cannot be idempotent.

\subsection{The nondecreasing case}

Starting from Proposition~\ref{prop:eoscm} we deduce immediately the following characterizations.

\begin{proposition}\label{prop:eondscm}
Assume that $E$ is open. Then $F\colon E^n\to \R$ is a nondecreasing comparison meaningful function on independent ordinal scales if and only if there
exist $k\in [n]$ and a strictly increasing or constant function $g\colon E\to\R$ such that $F=g\circ {\rm P}_k$.
\end{proposition}

When $E$ is not open, we have the following (see \cite[Corollary~4.2]{MarMesRuc05}).

\begin{theorem}
$F\colon E^n\to \R$ is a nondecreasing comparison meaningful function on independent ordinal scales if and only if there exists a mapping
$\xi\colon \isen\to [n]$ and a nondecreasing mapping $\gamma\colon\isen\to \mathcal{G}(E)$ such that $F|_I=(\gamma(I)\circ {\rm P}_{\xi(I)})|_I$ for all
$I\in\isen$, where if $\gamma(I)=\gamma(I')$ then also $\xi(I)=\xi(I')$ (unless $\gamma(I)=\gamma(I')$ is constant).
\end{theorem}

\subsection{The continuous case}

As we already mentioned in \S\ref{sec:sdsgf}, comparison meaningfulness on independent ordinal scales reduces any lattice polynomial function to
a projection function onto one coordinate. From this result we deduce immediately the following characterizations; see \cite[\S{5}]{Mar02c}.

\begin{theorem}\label{thm:coscm}
$F\colon E^n\to \R$ is a continuous and comparison meaningful function on independent ordinal scales if and only if there exist $k\in [n]$ and a
continuous and strictly monotonic or constant function $g\colon E\to\R$ such that $F=g\circ {\rm P}_k$.
\end{theorem}

\begin{corollary}
$F\colon E^n\to \R$ is a continuous, idempotent, and comparison meaningful function on independent ordinal scales if and only if it is a projection
function.
\end{corollary}

\begin{remark}
The result in Theorem~\ref{thm:coscm} was proved first by Kim~\cite[Corollary~1.2]{Kim90} in $\R^n$ (see Osborne~\cite{Osb70} for preliminary
results).
\end{remark}

It follows from Theorem~\ref{thm:coscm} that, if $n\geqslant 2$, any symmetric, continuous, and comparison meaningful function on independent
ordinal scales is necessarily a constant function.

%---------------------------------------------------------------------------------------------- Section 7
\section{Aggregation on finite chains by chain independent functions}
\label{sec:aod}

In this final section, mainly based on a paper by the authors~\cite{MarMes04}, we give interpretations of order invariance and comparison
meaningfulness properties in the setting of aggregation on finite chains (i.e., totally ordered finite sets). These interpretations show that
the order invariant functions and comparison meaningful functions always have isomorphic discrete representatives defined on finite chains.
These discrete functions do not depend on the chains on which they are defined.

\subsection{Introduction}

Let $A$ be a set of \emph{alternatives} (objects, individuals, etc.) and consider an open real interval $E$, possibly
unbounded.\footnote{Without loss of generality, we can assume that $E=(0,1)$ or $E=\R$.} In representational measurement
theory~\cite{Rob79,Rob94}, a \emph{scale of measurement} can be seen as a mapping $h\colon A\to E$ that assigns a real number to each element of $A$
according to some attribute or criterion.\footnote{A criterion is an attribute defined in a preference-ordered domain.} As already mentioned in
the introduction, such a scale is an ordinal scale if any other acceptable version of it is of the form $\phi\circ h$ for some strictly
increasing function $\phi\colon E\to E$.

An ordinal scale is finite if $\ran(h)$ is a finite subset of $E$, that is of the form
$$
\ran(h)=\{e_1 < e_2 < \cdots < e_k\},
$$
where the values $e_1,e_2,\ldots,e_k$ represent the possible rating benchmarks defined along some ordinal criterion. We shall assume throughout
that $|\ran(h)|=k\geqslant 2$.

Since the values $e_1,e_2,\ldots,e_k$ of the scale are defined up to order, that is, within a strictly increasing function $\phi\colon E\to E$, we can
simply replace $\ran(h)$ with a finite chain $(S,\preccurlyeq)$ of $k$ elements, that is,
$$
S=\{s_1 \prec s_2 \prec \cdots \prec s_k\},
$$
where $\preccurlyeq$ represents a total order on $S$ and $\prec$ represents its asymmetric part. In this representation we denote by $s_*:=s_1$
(resp.\ $s^*:=s_k$) the bottom element (resp.\ top element) of the chain.

\begin{example}
Consider the problem of evaluating a commodity by a consumer according to a given ordinal criterion. Typically this evaluation is done by rating
the product on a finite ordinal scale. For instance we could consider the following rating benchmarks:
$$
1=\mathrm{Bad},~2=\mathrm{Weak},~3=\mathrm{Fair},~4=\mathrm{Good},~5=\mathrm{Excellent}.
$$
Since the scale values are determined only up to order, this scale can be replaced with a finite chain $S=\{\mathrm{B} \prec \mathrm{W} \prec \mathrm{F} \prec \mathrm{G} \prec \mathrm{E}\}$ whose elements B, W, F, G, E refer to the following linguistic terms: {\it bad}, {\it weak}, {\it fair}, {\it good}, {\it excellent}.
\end{example}

It is well known (see \cite[Chapter~1]{KraLucSupTve71}) that the total order $\preccurlyeq$ defined on $S$ can always be
numerically represented in $E$ by means of an isomorphism (a strictly increasing function) $f\colon S\to E$ such that
$$
s_i\preccurlyeq s_j \quad\Leftrightarrow \quad f(s_i)\leqslant f(s_j) \qquad (s_i,s_j\in S).
$$
Moreover, just as for the mapping $h$, the isomorphism $f$ is defined up to a strictly increasing function $\phi\colon E\to E$. That is, with $f$ all
isomorphisms $f'=\phi\circ f$ (and only these) represent the same order on $S$.

By choosing $f$ so that $f(s_i)=e_i$ for $i=1,\ldots,k$, we immediately see that the elements of $A$ can be ordinally evaluated not only by
means of the numerical mapping $h\colon A\to \ran(h)$ but also by the non-numerical mapping $h_S\colon A\to S$, defined by $h_S:=f^{-1}\circ h$. The
following diagram illustrates the relationship among the mappings, where $h$ and $f$ are defined within a strictly increasing function
$\phi\colon E\to E$:
$$
\xymatrix{
    A \ar[r]^{\phi\circ h} \ar[rd]_{h_S} & E  \\
     & S \ar[u]_{\phi\circ f}
}
$$

We may also consider non-open intervals $E$ with the natural condition that if $f(s)=\inf E\in E$ for some $s\in S$ then $s=s_*$, and similarly,
if $f(s)=\sup E\in E$ for some $s\in S$ then $s=s^*$. In that case, the isomorphism $f$ is required to be \emph{endpoint preserving}, that is,
if $\inf E\in E$ (resp.\ $\sup E\in E$) then $f(s_*)=\inf E$ (resp.\ $f(s^*)=\sup E$), regardless of the chain $(S,\preccurlyeq)$
considered.\footnote{This amounts to assuming that all the chains considered have a common bottom element $s_*$ (resp.\ a common top element
$s^*$) whose numerical representation is $\inf E$ (resp.\ $\sup E$).} Consequently also all the functions $\phi\colon E\to E$ must be endpoint
preserving in the sense that $\phi(x)=x$ for all $x\in B[E]$. Due to the finiteness of the ordinal scales, we may even assume that the functions
$\phi$ are continuous, which amounts to assuming that they all belong to $\Phi[E]$.

The endpoint preservation assumption of $f$ (and hence of $\phi$) clarifies why we consider numerical representations in an interval $E$ of
$\mathbb{R}$, possibly non-open, rather than $\mathbb{R}$ itself.\footnote{If $E$ is closed, one typically chooses $E=[0,1]$ or
$E=\R\cup\{-\infty,\infty\}$.}

In this section, the set of all endpoint preserving isomorphisms $f\colon S\to E$ is denoted $\mathrm{F}[S,E]$. The diagonal restriction of
$\mathrm{F}[S,E]^n$ is the set
$$
\mathrm{F}_n[S,E]:=\{\underbrace{(f,\ldots,f)}_{n} : f\in\mathrm{F}[S,E]\}.
$$
Finally, for any $\ba\in S^n$ and any $\bff\in\mathrm{F}[S,E]^n$, the symbol $\bff(\ba)$ denotes the vector
$\big(f_1(a_1),\ldots,f_n(a_n)\big)$.

\subsection{Aggregation by order invariant functions}
\label{sec:aif}

Suppose we have $n$ evaluations expressed in a finite chain $(S,\preccurlyeq)$, with $|S|=k\geqslant 2$. To aggregate these evaluations and
obtain an overall evaluation in the same chain, we can use a discrete aggregation function $G\colon S^n\to S$, which is a ranking function sorting
$k^n$ $n$-tuples into $k$ classes. (Here, ``discrete" means that the domain of the function $G$ is a discrete set.)

Among all the possible aggregation functions, we could choose one that is ``independent'' of the chain used.\footnote{For example, we could use
any lattice polynomial function, which does not depend on the chain used.} Such a \emph{chain independent}\/ aggregation function is necessarily
based on a numerical function $F\colon E^n\to E$ that can be represented in any finite chain $(S,\preccurlyeq)$ by a discrete analog $G\colon S^n\to S$ in
the sense that the following identity
$$
F(x_1,\ldots,x_n)=f\big(G\big(f^{-1}(x_1),\ldots,f^{-1}(x_n)\big)\big)\qquad (\bx\in E^n)
$$
holds for all isomorphisms $f\in\mathrm{F}[S,E]$.

As the following theorem shows~\cite[Proposition~4.1]{MarMes04}, this condition completely characterizes the order invariant functions.

\begin{theorem}\label{thm:si}
$F\colon E^n\to E$ is an order invariant function if and only if, for any finite chain $(S,\preccurlyeq)$, there exists an aggregation function
$G\colon S^n\to S$ such that, for any $f\in\mathrm{F}[S,E]$, we have
\begin{equation}\label{eq:si}
F\big(f(a_1),\ldots,f(a_n)\big)=f\big(G(a_1,\ldots,a_n)\big) \qquad (\ba\in S^n).
\end{equation}
\end{theorem}

Thus, an order invariant function is characterized by the fact that it can always be represented by a discrete aggregation function $G\colon S^n\to S$
on any finite chain $(S,\preccurlyeq)$, regardless of the cardinality of this chain.\footnote{It is important to remember that considering a
discrete function $G\colon S^n\to S$, where $(S,\preccurlyeq)$ is a given chain, is not equivalent to considering an order invariant function
$F\colon E^n\to E$. Indeed, defining an order invariant function is much more restrictive since such a function should be independent of any scale.
For instance, if $n=2$ and $E$ is open, we see by Proposition~\ref{prop:eoin} that there are only four order invariant functions (namely $x_1$,
$x_2$, $x_1\wedge x_2$, and $x_1\vee x_2$) while the number of possible discrete functions $G\colon S^2\to S$ is clearly $k^{k^2}$, where $k=|S|$.} It
is informative to represent the identity (\ref{eq:si}) by the following commutative diagram, where $\bff:=(f,\ldots,f)$:
$$
\xymatrix{
    E^n \ar[r]^F & E \\
    S^n \ar[r]_G \ar[u]^{\bff} & S \ar[u]_{f}
}
$$

It is clear from (\ref{eq:si}) that the discrete function $G$ representing $F$ in $(S,\preccurlyeq)$ is uniquely determined and, in some sense,
is isomorphic to the ``restriction'' of $F$ to $S^n$. For example, if $n=2$ and $F(\bx)=x_1\wedge x_2$ (resp.\ $F(\bx)=\inf E$) then the unique
representative $G$ of $F$ is defined by $G(\ba)=a_1\wedge a_2$ (resp.\ $G(\mathbf{a})=s_*$).

Evidently an order invariant function is nondecreasing if and only if its discrete representative is nondecreasing. Another property that might
be required on order invariant functions is continuity (see \S\ref{sec:OIcont}) whose discrete counterpart, called \emph{smoothness}, is defined
as follows (see \cite{GodSie88}).\footnote{Fodor~\cite[Theorem~2]{Fod00} (see \cite{MarMes04} for the
general case) showed that the smoothness condition is equivalent to the discrete version of the intermediate value
theorem~\cite[Lemma~1]{FunFu75}.}

\begin{definition}\label{de:smooth}
Consider $(n+1)$ finite chains $(S_0,\preccurlyeq_{S_0}),\ldots,(S_n,\preccurlyeq_{S_n})$. A discrete function $G\colon\bigtimes_{i=1}^n S_i\to S_0$
is said to be \emph{smooth} if, for any $\ba,\bb\in \bigtimes_{i=1}^n S_i$, the elements $G(\ba)$ and $G(\bb)$ are equal or neighboring whenever
there exists $j\in [n]$ such that $a_j$ and $b_j$ are neighboring and $a_i=b_i$ for all $i\neq j$.
\end{definition}

The following important result~\cite[Proposition~5.1]{MarMes04} relates the continuity property of order invariant functions to the smoothness
condition of its discrete representatives, thus making continuity sensible and even appealing for order invariant functions.

\begin{proposition}\label{prop:cisr}
An order invariant function $F\colon E^n\to E$ is continuous if and only if it is represented only by smooth discrete aggregation functions.
\end{proposition}

\subsection{Aggregation by comparison meaningful functions on a single ordinal scale}
\label{sec:aif2}

Consider the more general situation where the evaluations to be aggregated are expressed in the same finite chain $(S,\preccurlyeq_S)$ and the
overall evaluation is expressed in a finite chain $(T,\preccurlyeq_T)$, possibly different from $(S,\preccurlyeq_S)$. Again, we can consider
aggregation functions $G\colon S^n\to T$ and, among them, we might want to choose aggregation functions that are independent of the chains used.

As the following theorem shows~\cite[Proposition~4.3]{MarMes04}, such chain independent functions are constructed from numerical functions
$F\colon E^n\to\R$ that are exactly the comparison meaningful functions on a single ordinal scale.

\begin{theorem}\label{thm:si2}
$F\colon E^n\to \R$ is a comparison meaningful function on a single ordinal scale if and only if, for any finite chain $(S,\preccurlyeq_S)$, there
exists a finite chain $(T,\preccurlyeq_T)$ and a surjective aggregation function $G\colon S^n\to T$ such that, for any $\bff\in\mathrm{F}_n[S,E]$,
there is an isomorphism $g_{\bff}\colon T\to\R$ such that
\begin{equation}\label{eq:iuif}
F\big(\bff(\ba)\big)=g_{\bff}\big(G(\ba)\big) \qquad (\ba\in S^n).
\end{equation}
\end{theorem}

Thus, a comparison meaningful function on a single ordinal scale is characterized by the fact that it can always be represented by a discrete
aggregation function $G\colon S^n\to T$ on any finite chain $(S,\preccurlyeq)$, regardless of the cardinality of this chain. The identity
(\ref{eq:iuif}) can be graphically represented by the following commutative diagram
$$
\xymatrix{
    E^n \ar[r]^F & \mathbb{R} \\
    S^n \ar[r]_G \ar[u]^{\bff} & T \ar[u]_{g_{\bff}}
}
$$

It can be easily shown~\cite[\S{4.2}]{MarMes04} that, given a comparison meaningful function $F\colon E^n\to \R$ on a single ordinal scale and a
finite chain $(S,\preccurlyeq_S)$, the output chain $(T,\preccurlyeq_T)$ and the functions $G\colon S^n\to T$ and $g_{\bff}\colon T\to\R$ are uniquely
determined.

The analog of Proposition~\ref{prop:cisr} can be stated as follows~\cite[Proposition~5.2]{MarMes04}. Unfortunately here we no longer have a
necessary and sufficient condition.

\begin{proposition}
A continuous comparison meaningful function $F\colon E^n\to \R$ on a single ordinal scale is represented only by smooth discrete aggregation
functions.
\end{proposition}

\subsection{Aggregation by comparison meaningful functions on independent ordinal scales}
\label{sec:aif3}

We now assume that the $n$ evaluations are expressed in independent finite chains $(S_i,\preccurlyeq_{S_i})$, $i=1,\ldots,n$, and that the
overall evaluation is expressed in a finite chain $(T,\preccurlyeq_T)$. We can consider aggregation functions $G\colon \bigtimes_{i=1}^n S_i\to T$
and, among them, we might want to choose aggregation functions that are independent of the chains used.

As the following theorem shows~\cite[Proposition~4.6]{MarMes04}, such chain independent functions are constructed from numerical functions
$F\colon E^n\to\R$ that are exactly the comparison meaningful functions on independent ordinal scales.

\begin{theorem}\label{thm:si3}
$F\colon E^n\to \R$ is a comparison meaningful function on independent ordinal scales if and only if, for any finite chains
$(S_i,\preccurlyeq_{S_i})$, $i=1,\ldots,n$, there exists a finite chain $(T,\preccurlyeq_T)$ and a surjective aggregation function
$G\colon \bigtimes_{i=1}^n S_i\to T$ such that, for any $\bff\in\mathrm{F}[S,E]^n$, there is a isomorphism $g_{\bff}\colon T\to\R$ such that
$$%\begin{equation}\label{eq:iuif3}
F\big(\bff(\ba)\big)=g_{\bff}\big(G(\ba)\big) \qquad (\ba\in \bigtimes_{i=1}^n S_i).
$$%\end{equation}
\end{theorem}

Thus, a comparison meaningful function on independent ordinal scales is characterized by the fact that it can always be represented by a
discrete aggregation function $G\colon\bigtimes_{i=1}^n S_i\to T$, regardless of the cardinality of the chains considered. Here the commutative
diagram is given by
$$
\xymatrix{
    E^n \ar[r]^F & \mathbb{R} \\
    \bigtimes\limits_{i=1}^n S_i \ar[r]_G \ar[u]^{\bff} & T \ar[u]_{g_{\bff}}
}
$$

Here again, it can be easily shown~\cite[\S{4.3}]{MarMes04} that, given a comparison meaningful function $F\colon E^n\to \R$ on independent ordinal
scales and $n$ finite chains $(S_i,\preccurlyeq_{S_i})$, $i=1,\ldots,n$, the output chain $(T,\preccurlyeq_T)$ and the functions
$G\colon\bigtimes_{i=1}^n S_i\to T$ and $g_{\bff}\colon T\to\R$ are uniquely determined.

Regarding continuous comparison meaningful functions, we have the following result~\cite[Proposition~5.3]{MarMes04}.

\begin{proposition}
A continuous comparison meaningful function $F\colon E^n\to \R$ on independent ordinal scales is represented only by smooth discrete aggregation
functions.
\end{proposition}

\section*{Acknowledgments}

Radko Mesiar greatly acknowledges the support of grants APVV-0375-06 and VEGA 1/4209/07.

%\bibliographystyle{abbrv}   % styles: plain, unsrt, alpha, abbrv, ieeetr, acm, siam, apalike, amsplain,...
%\bibliography{ReferencesMarichal}

\begin{thebibliography}{10}

\bibitem{Acz87}
J.~Acz{\'e}l.
\newblock {\em A short course on functional equations}.
\newblock Theory and Decision Library. Series B: Mathematical and Statistical
  Methods. D. Reidel Publishing Co., Dordrecht, 1987.

\bibitem{AczGroSch94}
J.~Acz{\'e}l, D.~Gronau, and J.~Schwaiger.
\newblock Increasing solutions of the homogeneity equation and of similar
  equations.
\newblock {\em J. Math. Anal. Appl.}, 182(2):436--464, 1994.

\bibitem{AczRob89}
J.~Acz{\'e}l and F.~S. Roberts.
\newblock On the possible merging functions.
\newblock {\em Math. Social Sci.}, 17(3):205--243, 1989.

\bibitem{AczRobRos86}
J.~Acz{\'e}l, F.~S. Roberts, and Z.~Rosenbaum.
\newblock On scientific laws without dimensional constants.
\newblock {\em J. Math. Anal. Appl.}, 119(1-2):389--416, 1986.

\bibitem{BarDre04}
L.~Bart{\l}omiejczyk and J.~Drewniak.
\newblock A characterization of sets and operations invariant under bijections.
\newblock {\em Aequationes Math.}, 68(1):1--9, 2004.

\bibitem{Bir67}
G.~Birkhoff.
\newblock {\em Lattice theory}.
\newblock Third edition. American Mathematical Society Colloquium Publications,
  Vol. XXV. American Mathematical Society, Providence, R.I., 1967.

\bibitem{BouPir97}
D.~Bouyssou and M.~Pirlot.
\newblock Choosing and ranking on the basis of fuzzy preference relations with
  the ``min in favor''.
\newblock In {\em Multiple criteria decision making (Hagen, 1995)}, volume 448
  of {\em Lecture Notes in Econom. and Math. Systems}, pages 115--127.
  Springer, Berlin, 1997.

\bibitem{Cau21}
A.~L. Cauchy.
\newblock {\em Cours d'analyse de l'Ecole Royale Polytechnique, Vol. I. Analyse
  alg\'ebrique}.
\newblock Debure, Paris, 1821.

\bibitem{DeBMes03}
B.~De~Baets and R.~Mesiar.
\newblock Discrete triangular norms.
\newblock In S.~E. Rodabaugh and E.~P. Klement, editors, {\em Topological and
  Algebraic Structures in Fuzzy Sets}, pages 389--400. Kluwer Academic
  Publishers, 2003.

\bibitem{FodRou95}
J.~Fodor and M.~Roubens.
\newblock {On meaningfulness of means}.
\newblock {\em J. Comput. Appl. Math.}, 64(1-2):103--115, 1995.

\bibitem{Fod00}
J.~C. Fodor.
\newblock Smooth associative operations on finite ordinal scales.
\newblock {\em IEEE Trans. Fuzzy Syst.}, 8(6):791--795, 2000.

\bibitem{FunFu75}
L.~W. Fung and K.~S. Fu.
\newblock An axiomatic approach to rational decision making in a fuzzy
  environment.
\newblock In {\em Fuzzy sets and their applications to cognitive and decision
  processes (Proc. U. S.-Japan Sem., Univ. Calif., Berkeley, Calif., 1974)},
  pages 227--256. Academic Press, New York, 1975.

\bibitem{GarMar08}
J.~L. Garc\'{\i}a-Lapresta and R.~A. Marques~Pereira.
\newblock The self-dual core and the anti-self-dual remainder of an aggregation
  operator.
\newblock {\em Fuzzy Sets and Systems}, 159(1):47--62, 2008.

\bibitem{GodSie88}
L.~Godo and C.~Sierra.
\newblock A new approach to connective generation in the framework of expert
  systems using fuzzy logic.
\newblock In {\em Proc. 18th Int. Symposium on Multiple-Valued Logic}, pages
  157--162, Palma de Mallorca, Spain, 24-26 May 1988.

\bibitem{Grae03}
G.~Gr\"atzer.
\newblock {\em General lattice theory}.
\newblock Birkh\"auser Verlag, Berlin, 2003.
\newblock Second edition.

\bibitem{Gri07}
P.~A. Grillet.
\newblock {\em Abstract algebra}, volume 242 of {\em Graduate Texts in
  Mathematics}.
\newblock Springer, New York, second edition, 2007.

\bibitem{HarLitPol52}
G.~Hardy, J.~Littlewood, and G.~P\'olya.
\newblock {\em {Inequalities. 2nd ed.}}
\newblock {Cambridge, Engl.: At the University Press}, 1952.

\bibitem{KamOvc95}
Y.~Kamen and S.~Ovchinnikov.
\newblock Meaningful means on ordered sets.
\newblock In B.~Bouchon-Meunier, R.~R. Yager, and L.~A. Zadeh, editors, {\em
  Fuzzy logic and soft computing}, pages 189--193. World Sci. Publishing, River
  Edge, NJ, 1995.

\bibitem{Kim90}
S.-R. Kim.
\newblock On the possible scientific laws.
\newblock {\em Math. Social Sci.}, 20(1):19--36, 1990.

\bibitem{Kle69}
D.~Kleitman.
\newblock {On {D}edekind's problem: the number of monotone {B}oolean
  functions}.
\newblock {\em Proc. Am. Math. Soc.}, 21:677--682, 1969.

\bibitem{KleMar75}
D.~Kleitman and G.~Markowsky.
\newblock {On {D}edekind's problem: The number of isotone {B}oolean functions.
  II.}
\newblock {\em Trans. Am. Math. Soc.}, 213:373--390, 1975.

\bibitem{KraLucSupTve71}
D.~H. Krantz, R.~D. Luce, P.~Suppes, and A.~Tversky.
\newblock {\em Foundations of measurement. Vol. I: Additive and polynomial
  representations}.
\newblock Academic Press, New York, 1971.

\bibitem{KucChoGer90}
M.~Kuczma, B.~Choczewski, and R.~Ger.
\newblock {\em Iterative Functional Equations}.
\newblock Cambridge University Press, Cambridge, UK, 1990.

\bibitem{Luc59}
R.~D. Luce.
\newblock On the possible psychophysical laws.
\newblock {\em Psych. Rev.}, 66:81--95, 1959.

\bibitem{LucKraSupTve90}
R.~D. Luce, D.~H. Krantz, P.~Suppes, and A.~Tversky.
\newblock {\em Foundations of measurement. Vol. III: Representation,
  axiomatization, and invariance}.
\newblock Academic Press Inc., San Diego, CA, 1990.

\bibitem{Mar98}
J.-L. Marichal.
\newblock {\em Aggregation operators for multicriteria decision aid}.
\newblock PhD thesis, Institute of Mathematics, University of Li\`ege, Li\`ege,
  Belgium, December 1998.

\bibitem{Mar01}
J.-L. Marichal.
\newblock An axiomatic approach of the discrete {S}ugeno integral as a tool to
  aggregate interacting criteria in a qualitative framework.
\newblock {\em IEEE Trans. Fuzzy Syst.}, 9(1):164--172, 2001.

\bibitem{Mar02c}
J.-L. Marichal.
\newblock On order invariant synthesizing functions.
\newblock {\em J. Math. Psych.}, 46(6):661--676, 2002.

\bibitem{MarMat01}
J.-L. Marichal and P.~Mathonet.
\newblock On comparison meaningfulness of aggregation functions.
\newblock {\em J. Math. Psych.}, 45(2):213--223, 2001.

\bibitem{MarMes04}
J.-L. Marichal and R.~Mesiar.
\newblock Aggregation on finite ordinal scales by scale independent functions.
\newblock {\em Order}, 21(2):155--180, 2004.

\bibitem{MarMesRuc05}
J.-L. Marichal, R.~Mesiar, and T.~R\"uckschlossov\'a.
\newblock A complete description of comparison meaningful functions.
\newblock {\em Aequationes Math.}, 69(3):309--320, 2005.

\bibitem{MarRou93}
J.-L. Marichal and M.~Roubens.
\newblock Characterization of some stable aggregation functions.
\newblock In {\em Proc.\ 1st Int.\ Conf.\ on Industrial Engineering and
  Production Management (IEPM'93)}, pages 187--196, Mons, Belgium, June 1993.

\bibitem{Mes01}
R.~Mesiar.
\newblock Scale invariant operators.
\newblock In {\em Proc. Int. Conf. in Fuzzy Logic and Technology
  (EUSFLAT'2001)}, pages 479--481, Leicester, September 5-7 2001.

\bibitem{MesRuc04}
R.~Mesiar and T.~R\"uckschlossov\'a.
\newblock Characterization of invariant aggregation operators.
\newblock {\em Fuzzy Sets and Systems}, 142(1):63--73, 2004.

\bibitem{Nar02}
L.~Narens.
\newblock {\em Theories of meaningfulness}.
\newblock Lawrence Erlbaum Associates, Mahwah, N.J., 2002.

\bibitem{Orl81}
A.~I. Orlov.
\newblock The connection between mean values and the admissible transformations
  of scale.
\newblock {\em Math. Notes}, 30:774--778, 1981.

\bibitem{Osb70}
D.~K. Osborne.
\newblock {Further extensions of a theorem of dimensional analysis}.
\newblock {\em J. Math. Psychol.}, 7:236--242, 1970.

\bibitem{Ovc96}
S.~Ovchinnikov.
\newblock Means on ordered sets.
\newblock {\em Math. Social Sci.}, 32(1):39--56, 1996.

\bibitem{Ovc98c}
S.~Ovchinnikov.
\newblock Invariant functions on simple orders.
\newblock {\em Order}, 14(4):365--371, 1998.

\bibitem{Ovc98}
S.~Ovchinnikov.
\newblock On ordinal {OWA} operators.
\newblock In {\em Proc.\ 7th Int.\ Conf.\ on Information Processing and
  Management of Uncertainty in Knowledge-Based Systems (IPMU'98)}, pages
  511--514, Paris, 1998.

\bibitem{OvcDuk00}
S.~Ovchinnikov and A.~Dukhovny.
\newblock {Integral representation of invariant functionals}.
\newblock {\em J. Math. Anal. Appl.}, 244(1):228--232, 2000.

\bibitem{OvcDuk02}
S.~Ovchinnikov and A.~Dukhovny.
\newblock {On order invariant aggregation functionals}.
\newblock {\em J. Math. Psychol.}, 46(1):12--18, 2002.

\bibitem{RobRos94}
F.~Roberts and Z.~Rosenbaum.
\newblock The meaningfulness of ordinal comparisons for general order
  relational systems.
\newblock In P.~Humphreys, editor, {\em Patrick Suppes: Scientific Philosopher,
  Vol. 2}, pages 251--274. Kluwer, The Netherlands, 1994.

\bibitem{Rob79}
F.~S. Roberts.
\newblock {\em Measurement theory, with applications to decision-making,
  utility and the social sciences}, volume~7 of {\em Encyclopedia of
  Mathematics and its Applications}.
\newblock Addison-Wesley Publishing Co., Reading, Mass., 1979.

\bibitem{Rob94}
F.~S. Roberts.
\newblock Limitations on conclusions using scales of measurement.
\newblock In S.~M. Pollock, M.~H. Rothkopf, and A.~Barnett, editors, {\em
  Operations Research and the Public Sector}, pages 621--671. Elsevier,
  Amsterdam, 1994.

\bibitem{Yan89}
E.~B. Yanovskaya.
\newblock Group choice rules in problems with interpersonal preference
  comparisons.
\newblock {\em Automat. Remote Control}, 50(6):822--830, 1989.

\end{thebibliography}

\end{document}